\documentclass[10pt,twoside]{article}
\usepackage{graphicx,subfigure}
\usepackage{amsfonts}
\usepackage{mathrsfs}
\usepackage{amssymb}
\usepackage{graphicx}
\usepackage{amsmath,amssymb}
\usepackage[mathscr]{eucal}
\begin{document}

\centerline{\Large\bf Uniform asymptotic stability of solutions of
fractional  }\vspace{0.3cm} \centerline {\Large\bf  functional
differential equations  \footnote{Research supported by the National
Natural Science Foundation of China under Grant No. 10801066, and
the Fundamental Research Funds for the Central Universities under
Grant No. lzujbky-2011-47 and  No. lzujbky-2012-k26. }}
 \vskip 3pt \centerline{{\bf Yajing Li}~~ and ~~{\bf Yejuan Wang\footnote{Corresponding
author.\\
 \indent E-mail addresses: wangyj@lzu.edu.cn (Yejuan Wang), liyajing11@st.lzu.edu.cn
(Yajing Li).}}} \vskip 3pt

\centerline{\footnotesize\it School of Mathematics and Statistics,
Lanzhou University, Lanzhou 730000,
 PR China}

\bigbreak \noindent{\bf Abstract} {In this paper, some global
existence and uniform asymptotic stability results for fractional
functional  differential equations are proved. It is worthy
mentioning that when $\alpha=1$ the initial value problem
 (1.1)   reduces to a classical
dissipative differential equation with delays in \cite{Caraballo2}.

}

\vskip 5pt
 \noindent{{\it Keywords:} Functional differential equation;
 Fractional derivative; Asymptotic stability; Global existence.}

\noindent{\section{Introduction}}

Consider the initial value problem (IVP for short) of the following
fractional functional  differential equation:
$$
\left\{\begin{array}{l} \displaystyle D^\alpha \left[y(t)e^{\beta
t}\right] =f(t,y_t)e^{\beta t},
\quad t\in[t_0, \infty), \, t_0 \geqslant 0, \, 0<\alpha <1,\\
y(t)=\phi(t),\quad t_0-h \leqslant t \leqslant t_0,
\end{array}
\right. \eqno(1.1)
$$
where $D^\alpha$ is the  Caputo fractional derivative, $\beta >0$,
$f : J \times C([-h, 0], \mathbb{R}) \rightarrow \mathbb{R}$, where
$J=[t_0, \infty)$, is a given function satisfying some assumptions
that will be specified later, $h >0$, and  $\phi \in C([t_0-h, t_0],
\mathbb{R})$. If $y \in C([t_0-h, \infty), \mathbb{R})$, then for
any $t \in [t_0, \infty)$, define $y_t$ by
$$y_t(\theta)=y(t+\theta), \quad  \theta \in [-h, 0].
$$

The study of retarded differential equations  is an important area
of applied mathematics due to physical reasons, non-instant
transmission phenomena, memory processes, and specially biological
motivations (see, e.g., \cite{Caraballo2, Hale1, Kuang, Mackey}).
Fractional differential equations have attracted much attention
recently (see, for example, \cite{Benchohra, Benchohra1, Diethelm,
Diethelm2, Kilbas, Messaoudi, Podlubny} and the references cited
therein  for the applications in various sciences such as physics,
mechanics, chemistry, engineering, etc).

Some attractivity  results for fractional functional differential
equations and nonlinear functional integral equations are obtained
by using the fixed point theory; see \cite{Chen1, Chen2, Dhage1,
Dhage2, Dhage3} and references therein. Global asymptotic stability
of solutions of a functional integral equation is discussed in
\cite{Bana}, however there is no work on uniform asymptotic
stability of solutions of fractional functional differential
equation.  It is our intention here to show the global existence and
uniform asymptotic stability of the  fractional functional
differential equation (1.1).

We organize the paper as follows. In Section 2, we recall some
necessary concepts and results. In Section 3 we give the global
existence and uniform asymptotic stability of fractional functional
 differential equations. Finally, two examples are given to illustrate our main results.

\bigbreak

\noindent{\section{Preliminaries}}

In this section, we introduce notations, definitions, and
preliminary facts which are used throughout this paper.

We consider  $BC := BC([t_0-h, \infty), \mathbb{R})$ the Banach
space of all bounded and continuous functions  from $[t_0-h,
\infty)$ into $\mathbb{R}$ with the norm
$$
\|y\|_{\infty} :=  \sup\left\{|y(t)| ~:~ t \in [t_0-h,
\infty)\right\}.
$$
 Let $\|y_t\|=\sup_{-h \leqslant \theta \leqslant 0} |y(t+\theta)|$
for $t \in J$.

Throughout this paper, we always assume that $f(t, x_t)$ satisfies
the following condition:
\begin{enumerate}
\item[$(H_0)$] $f(t, x_t)$ is Lebesgue measurable with respect to $t$
on $[t_0, \infty)$, and $f(t, \varphi)$ is continuous with respect
to $\varphi$ on $C([-h, 0], \mathbb{R})$.
\end{enumerate}
By condition $(H_0)$ and the technique used in \cite{Diethelm}, we
get the equivalent form of IVP (1.1) as:
$$
y(t)=\left\{\begin{array}{l} \displaystyle y(t_0)e^{-\beta
(t-t_0)}+\frac{1}{\Gamma(\alpha)}\int\limits_{t_0}^t
(t-s)^{\alpha-1}e^{-\beta (t-s)}f(s, y_s)ds,
\quad  t \geqslant t_0,\\
\phi(t), \quad t \in [t_0-h, t_0].
\end{array}
\right. \eqno(2.1)
$$

\vskip 5pt \noindent {\bf Definition 2.1.} {\sl We say that
solutions of IVP $(1.1)$  are uniformly asymptotically stable if for
any bounded subset $B$ of $C([-h, 0], \mathbb{R})$ and
$\varepsilon>0$, there exists a $T>0$ such that
$$
|y(t, t_0, \phi)-x(t, t_0, \psi)| \leqslant \varepsilon \quad \mbox{
for all } t \geqslant T \mbox{ and } \phi, \psi \in B.
$$

}

We recall the following generalization of Gronwall's lemma for
singular kernels \cite{Henry}, which will be used in the sequel.

\vskip 5pt \noindent {\bf Lemma 2.2.} {\sl Let $v: [t_0, b]
\rightarrow [0, +\infty)$ be a real function and $w(\cdot)$ is a
nonnegative, locally integrable function on $[t_0, b]$ and there are
constants $a > 0$ and $0 < \alpha < 1$ such that
$$
 v(t)\leqslant w(t)+a\int\limits_{t_0}^t\frac{v(s)}{(t-s)^\alpha}ds.
 $$
Then there exists a constant $K=K(\alpha)$ such that
$$
v(t)\leqslant w(t)+K
a\int\limits_{t_0}^t\frac{w(s)}{(t-s)^\alpha}ds,
$$
for every $t \in [t_0, b]$. }

\vskip 5pt \noindent {{\bf Theorem 2.3} (Leray-Schauder Fixed Point
Theorem).} {\sl  Let $P$ be a continuous and compact mapping of a
Banach space $X$ into itself, such that the set
$$
\{x \in X ~: ~ x=\lambda Px \mbox{ for some } 0 \leqslant \lambda
\leqslant 1\}
$$
is bounded. Then $P$  has a fixed point.
 }

\noindent{\section{FDEs of fractional order}}

In this section, we will investigate the IVP (1.1). Our first global
existence and uniform asymptotic stability result for the IVP (1.1)
is based on the Banach contradiction principle and Lemma 2.2.

\vskip 5pt \noindent {\bf Theorem 3.1.} {\sl
 Assume that $f(t, y_t)$ satisfies conditions $(H_0)$ and
\begin{enumerate}
\item[$(H_1)$] there exists  $l>0$ such that
$$
|f(t,u_t)-f(t,v_t)| \leqslant l\|u_t-v_t\| \eqno(3.1)
$$
 for  $t \in J$ and every  $u_t, v_t \in C([-h,
0], \mathbb{R})$. Moreover, the function $t \mapsto f(t, 0)$ is
bounded with $f_0=\sup_{t \geqslant t_0}|f(t, 0)|$.
\end{enumerate}
If
 $$l\left(\frac{(t_0+h)^{\alpha-1}e^{-\beta (t_0+h)}}{\beta\Gamma(\alpha)}+
\frac{(t_0+h)^{\alpha}}{\Gamma(\alpha+1)}\right)<1,\eqno(3.2)
$$
 then the  IVP $(1.1)$  has a unique solution in the space $BC$.
 Moreover, solutions of IVP $(1.1)$ are uniformly asymptotically stable.}

 \vskip 5pt \noindent {\bf Proof.} We divide the proof into two steps.

\noindent Step1. We  define  the operator $P : C([t_0-h, \infty),
\mathbb{R}) \rightarrow C([t_0-h, \infty), \mathbb{R})$
 by
$$
(Py)(t)=\left\{\begin{array}{l}\displaystyle y(t_0)e^{-\beta
(t-t_0)}+\frac{1}{\Gamma(\alpha)}\int\limits_{t_0}^t
(t-s)^{\alpha-1}e^{-\beta (t-s)}f(s,y_s)ds, \quad  t \geqslant t_0,\\
\phi(t),\quad  t\in[t_0-h, t_0].
\end{array}
\right. \eqno(3.3)
$$
The operator $P$ maps $BC$ into itself. Indeed for each $y\in BC$,
and for each $t \geqslant 2t_0+h$, it follows from $(H_1)$ that
$$
\begin{array}{l}
\displaystyle|(Py)(t)|  \leqslant \displaystyle |y(t_0)|e^{-\beta
(t-t_0)}+\frac{1}{\Gamma(\alpha)}\int\limits_{t_0}^t
(t-s)^{\alpha-1}e^{-\beta (t-s)}(f_0+l\|y_s\|)ds\\
~~ \leqslant  \displaystyle \|\phi\|e^{-\beta
(t-t_0)}+\frac{f_0+l\|y\|_{\infty}}{\Gamma(\alpha)}\left(\int\limits_{t_0}^{t-(t_0+h)}
(t_0+h)^{\alpha-1}e^{-\beta (t-s)}ds+\int\limits_{t-(t_0+h)}^t
(t-s)^{\alpha-1}ds\right)\\
 ~~ \leqslant  \displaystyle
\|\phi\|+(f_0+l\|y\|_{\infty})\left(\frac{(t_0+h)^{\alpha-1}e^{-\beta
(t_0+h)}}{\beta\Gamma(\alpha)}+
\frac{(t_0+h)^{\alpha}}{\Gamma(\alpha+1)}\right),
\end{array}
$$
for each $t \in [t_0, 2t_0+h]$ we have
$$
\displaystyle|(Py)(t)| \leqslant \|\phi\|+\frac{(f_0+l
\|y\|_{\infty})(t_0+h)^\alpha}{\Gamma(\alpha+1)},
$$
and consequently $P(y) \in BC$.

 Since $BC := BC([t_0-h, \infty),
\mathbb{R})$ is a Banach space with norm $\|\cdot\|_{\infty}$, we
shall show that $P : BC \rightarrow BC$ is a contraction map. Let
$y_1, y_2 \in BC$. Then we have for each $t \geqslant t_0$,
$$
\begin{array}{lll}
\displaystyle|(Py_1)(t)-(Py_2)(t)| &\leqslant&\displaystyle
\frac{1}{\Gamma(\alpha)}\int\limits_{t_0}^t
(t-s)^{\alpha-1}e^{-\beta (t-s)}|f(s,y_{1s})-f(s,y_{2s})|ds\\
\displaystyle & \leqslant
&\displaystyle\frac{l}{\Gamma(\alpha)}\int\limits_{t_0}^t
(t-s)^{\alpha-1}e^{-\beta (t-s)}\|y_{1s}-y_{2s}\|ds.
\end{array}\eqno(3.4)
$$
Therefore for any $t \geqslant 2t_0+h$,
$$
\begin{array}{l}
\displaystyle|(Py_1)(t)-(Py_2)(t)|\\
\displaystyle~~\leqslant\frac{l}{\Gamma(\alpha)}
\|y_1(\cdot)-y_2(\cdot)\|_{\infty}\left(\int\limits_{t_0}^{t-(t_0+h)}
(t_0+h)^{\alpha-1}e^{-\beta (t-s)}ds+\int\limits_{t-(t_0+h)}^t(t-s)^{\alpha-1}ds\right)\\
\displaystyle~~\leqslant l\left(\frac{(t_0+h)^{\alpha-1}e^{-\beta
(t_0+h)}}{\beta\Gamma(\alpha)}+
\frac{(t_0+h)^{\alpha}}{\Gamma(\alpha+1)}\right)\|y_1(\cdot)-y_2(\cdot)\|_{\infty},
\end{array}\eqno(3.5)
$$
and for $t_0-h \leqslant t \leqslant 2t_0+h$,
$$
\begin{array}{lll}
\displaystyle|(Py_1)(t)-(Py_2)(t)| & \leqslant & \displaystyle
\frac{l}{\Gamma(\alpha)}\|y_1(\cdot)-y_2(\cdot)\|_{\infty}\int\limits_{t_0}^t
(t-s)^{\alpha-1}ds\\
\displaystyle & \leqslant &\displaystyle
\frac{l(t_0+h)^\alpha}{\Gamma(\alpha+1)}\|y_1(\cdot)-y_2(\cdot)\|_{\infty},
\end{array}\eqno(3.6)
$$
and thus
$$
\begin{array}{l}
\displaystyle\|(Py_1)(\cdot)-(Py_2)(\cdot)\|_{\infty}\\
~~\displaystyle \leqslant l\left(\frac{(t_0+h)^{\alpha-1}e^{-\beta
(t_0+h)}}{\beta\Gamma(\alpha)}+
\frac{(t_0+h)^{\alpha}}{\Gamma(\alpha+1)}\right)\|y_1(\cdot)-y_2(\cdot)\|_{\infty}.
\end{array}
\eqno(3.7)
$$
Hence, (3.2) and (3.7) imply that the operator $P$ is a contraction.
Therefore, $P$ has a unique fixed point by Banach's contraction
principle.

\noindent Step2. For any two  solutions $x=x(t)$ and $y=y(t)$ of IVP
(1.1) corresponding to initial values $\psi$ and $\phi$,  by (2.1)
we can deduce that for all $t \geqslant t_0+h$ and all $\theta \in
[-h, 0]$,
$$
\begin{array}{lll}
\displaystyle |x(t+\theta)-y(t+\theta)|\leqslant\displaystyle
|x(t_0)-y(t_0)|e^{-\beta (t+\theta-t_0)}\\
\quad\quad\displaystyle+\frac{1}{\Gamma(\alpha)}\int\limits_{t_0}^{t+\theta}(t+\theta-s)^{\alpha-1}e^{-\beta (t+\theta-s)}|f(s,x_s)-f(s,y_s)|ds\\
~~\displaystyle \leqslant|x(t_0)-y(t_0)|e^{-\beta
(t+\theta-t_0)}+\frac{l}{\Gamma(\alpha)}
\int\limits_{t_0}^{t+\theta}(t+\theta-s)^{\alpha-1}e^{-\beta
(t+\theta-s)}\|x_s-y_s\|ds.
\end{array}\eqno(3.8)
$$
Then, it follows that
$$
\begin{array}{lll}
\displaystyle e^{\beta t}\|x_t-y_t\|\leqslant \displaystyle
|x(t_0)-y(t_0)|e^{\beta (h+t_0)} +\frac{le^{\beta
h}}{\Gamma(\alpha)}\int\limits_{t_0}^t(t-s)^{\alpha-1}e^{\beta
s}\|x_s-y_s\|ds.
\end{array}\eqno(3.9)
$$
Let $w(t)=e^{\beta t}\|x_t-y_t\|$. Then we have
$$
\begin{array}{lll}
\displaystyle w(t)\leqslant|x(t_0)-y(t_0)|e^{\beta (h+t_0)}
+\frac{le^{\beta h}}{\Gamma(\alpha)}\int\limits_{t_0}^t
(t-s)^{\alpha-1}w(s)ds.
\end{array}\eqno(3.10)
$$
Applying Lemma 2.2, one can see that there exists a constant $K$
such that
$$
\begin{array}{lll}
\displaystyle w(t)& \leqslant &\displaystyle|x(t_0)-y(t_0)|e^{\beta
(h+t_0)} +\frac{Kle^{\beta h}}{\Gamma(\alpha)}\int\limits_{t_0}^t
(t-s)^{\alpha-1}|x(t_0)-y(t_0)|e^{\beta (h+t_0)}ds\\
\displaystyle &\leqslant &\displaystyle|x(t_0)-y(t_0)|e^{\beta
(h+t_0)}\left(1+\frac{Kle^{\beta
h}}{\Gamma(\alpha+1)}(t-t_0)^\alpha\right).
\end{array}\eqno(3.11)
$$
Hence we obtain
$$
e^{\beta t}\|x_t-y_t\|=w(t)\leqslant|x(t_0)-y(t_0)|e^{\beta
(h+t_0)}\left(1+\frac{Kle^{\beta
h}}{\Gamma(\alpha+1)}(t-t_0)^\alpha\right),
$$
and thus for all $t \geqslant t_0+h$,
$$
\begin{array}{lll}
\displaystyle |x(t)-y(t)|\leqslant |x(t_0)-y(t_0)|e^{-\beta
(t-h-t_0)}\left(1+\frac{Kle^{\beta
h}}{\Gamma(\alpha+1)}(t-t_0)^\alpha\right),
\end{array}
$$
which implies that the solutions of IVP (1.1) are uniformly
asymptotically stable. $\Box$

Now we give global existence and uniform asymptotic stability
results based on the nonlinear alternative of Leray-Schauder type.

\vskip 5pt \noindent {\bf Theorem 3.2.} {\sl Assume that the
following hypotheses hold:
\begin{enumerate}
\item[$(H_2)$] $f$ is a continuous function;

\item[$(H_3)$] there exist positive functions  $k_1$, $k_2 \in BC([t_0, \infty), \mathbb{R}_+)$ such that
$$
|f(t, u_t)| \leqslant k_1(t)+k_2(t)\|u_t\|
$$
for $t \in J$ and every $u_t \in C([-h, 0], \mathbb{R})$;

\item[$(H_4)$] moreover,
assume that $K_1=\sup_{t \geqslant t_0}k_1(t)$,
$K_2=\sup_{t\geqslant t_0}k_2(t)$,
$$
 \lim \limits_{t \rightarrow \infty}\int\limits_{t_0}^t (t-s)^{\alpha-1}e^{-\beta
(t-s)}k_1(s)ds=0,
$$
and
$$
 \lim \limits_{t \rightarrow \infty}\int\limits_{t_0}^t (t-s)^{\alpha-1}e^{-\beta
(t-s)}k_2(s)ds=0.
$$

\end{enumerate}
Then the IVP $(1.1)$  admits a solution in the space $BC$.
 Moreover, solutions of IVP $(1.1)$ are uniformly asymptotically stable. }

\vskip 5pt \noindent {\bf Proof.} Let $P : C([t_0-h, \infty),
\mathbb{R}) \rightarrow C([t_0-h, \infty), \mathbb{R})$ be defined
as in (3.3). First we show that $P$ maps $BC$ into itself. Indeed,
the map  $P(y)$ is continuous on $[t_0-h,+\infty)$ for each $y\in
BC$, and for each $t \geqslant 2t_0+h$, $(H_2)$ implies that
$$
\begin{array}{l}
\displaystyle|(Py)(t)|  \leqslant \displaystyle |y(t_0)|e^{-\beta
(t-t_0)}+\frac{1}{\Gamma(\alpha)}\int\limits_{t_0}^t
(t-s)^{\alpha-1}e^{-\beta (t-s)}(K_1+K_2\|y_s\|)ds\\
~~ \leqslant  \displaystyle \|\phi\|e^{-\beta
(t-t_0)}+\frac{K_1+K_2\|y\|_{\infty}}{\Gamma(\alpha)}\left(\int\limits_{t_0}^{t-(t_0+h)}
(t_0+h)^{\alpha-1}e^{-\beta (t-s)}ds+\int\limits_{t-(t_0+h)}^t
(t-s)^{\alpha-1}ds\right)\\
 ~~ \leqslant  \displaystyle
\|\phi\|+(K_1+K_2\|y\|_{\infty})\left(\frac{(t_0+h)^{\alpha-1}e^{-\beta
(t_0+h)}}{\beta\Gamma(\alpha)}+
\frac{(t_0+h)^{\alpha}}{\Gamma(\alpha+1)}\right),
\end{array}\eqno(3.12)
$$
for each $t \in [t_0, 2t_0+h]$ we have
$$
\displaystyle|(Py)(t)| \leqslant \|\phi\|+\frac{(K_1+K_2
\|y\|_{\infty})(t_0+h)^\alpha}{\Gamma(\alpha+1)},\eqno(3.13)
$$
and for any $t \in [t_0-h, t_0]$,
$$
\displaystyle|(Py)(t)|\leqslant \|\phi\|.
$$
Thus,
$$
\|P(y)\|_{\infty} \leqslant
\|\phi\|+(K_1+K_2\|y\|_{\infty})\left(\frac{(t_0+h)^{\alpha-1}e^{-\beta
(t_0+h)}}{\beta\Gamma(\alpha)}+
\frac{(t_0+h)^{\alpha}}{\Gamma(\alpha+1)}\right),
$$
and consequently $P(y) \in BC$.

Next, we show that the operator $P$ is continuous and completely
continuous, and there exists an open set $U \subset BC$ with $y \neq
\lambda P(y)$ for $\lambda \in (0, 1)$ and $y \in \partial U$.

 \noindent Step 1. $P$ is continuous.

Let $\{y_n\}$ be a sequence such that $y_n \rightarrow y$ in $BC$.
Then   there exist $R>0$ and $N>0$ such that
$$
\|y_n\|_{\infty}+\|y\|_{\infty} < R, ~~~\forall n \geqslant N.
\eqno(3.14)
$$
Let $\varepsilon>0$ be given. Since $(H_4)$ holds, there is a real
number $T>0$ such that
$$
\frac{2}{\Gamma(\alpha)}\int\limits_{t_0}^t
(t-s)^{\alpha-1}e^{-\beta (t-s)}(k_1(s)+k_2(s)R)ds< \varepsilon
\eqno(3.15)
$$
for all $t \geqslant T$. Now we consider the following two cases.

\noindent Case 1: if $t \geqslant T$, then it follows from $(H_3)$
and  (3.14)-(3.15) that for $n$ sufficiently large
$$
\begin{array}{l}
\displaystyle|Py_n(t)-Py(t)|   \leqslant |y_n(t_0)-y(t_0)|e^{-\beta
(t-t_0)}+ \frac{1}{\Gamma(\alpha)}\int\limits_{t_0}^t
(t-s)^{\alpha-1}e^{-\beta (t-s)}|f(s,y_{ns})-f(s,y_{s})|ds\\
~~ \leqslant\displaystyle |y_n(t_0)-y(t_0)|+
\frac{2}{\Gamma(\alpha)}\int\limits_{t_0}^t
(t-s)^{\alpha-1}e^{-\beta (t-s)}(k_1(s)+k_2(s)R)ds<2\varepsilon.\\
\end{array}\eqno(3.16)
$$

\noindent Case 2: if $t_0 \leqslant t \leqslant T$,  since $f$ is a
continuous function, one has
$$
\begin{array}{lll}
\displaystyle |Py_n(t)-Py(t)|& \leqslant & \displaystyle
|y_n(t_0)-y(t_0)|+ \frac{1}{\Gamma(\alpha)}\int\limits_{t_0}^t
(t-s)^{\alpha-1}e^{-\beta (t-s)}|f(s,y_{ns})-f(s,y_{s})|ds\\
& \leqslant & \displaystyle |y_n(t_0)-y(t_0)|+
\frac{(T-t_0)^\alpha}{\Gamma(\alpha+1)}\sup \limits_{s \in [t_0,
T]}|f(s,y_{ns})-f(s,y_{s})|.
\end{array}\eqno(3.17)
$$
Note that $y_n \rightarrow y$ in $BC$. Hence (3.16) and (3.17) imply
that
$$
\|P(y_n)-P(y)\|_{\infty}\rightarrow 0 ~~~\mbox{ as } n \rightarrow
\infty.
$$

\noindent Step 2. $P$ maps bounded sets into bounded sets in $BC$.

Indeed, it is enough to show that for any $\eta>0$, there exists a
positive constant $\ell$ such that for  each $y \in B_\eta=\{y \in
BC : ||y||_{\infty}\leqslant \eta\}$ one has
$\|P(y)\|_{\infty}\leqslant \ell$. Let $y \in B_\eta$. Then we have
for each $t \geqslant 2t_0+h$,
$$
\begin{array}{l}
\displaystyle|(Py)(t)|\leqslant|y(t_0)|e^{-\beta (t-t_0)}+
\frac{1}{\Gamma(\alpha)}\int\limits_{t_0}^t
(t-s)^{\alpha-1}e^{-\beta (t-s)}|f(s,y_{s})|ds\\
\displaystyle~~\leqslant \eta
+\frac{K_1+K_2\|y\|_{\infty}}{\Gamma(\alpha)}\int\limits_{t_0}^t
(t-s)^{\alpha-1}e^{-\beta (t-s)}ds\\
~~\displaystyle\leqslant \eta+(K_1+K_2\eta)\left(
\frac{(t_0+h)^{\alpha-1}e^{-\beta(t_0+h)}}{\beta\Gamma(\alpha)}
+\frac{(t_0+h)^\alpha}{\Gamma(\alpha+1)}\right)=: \ell,
\end{array}
$$
and for each $t$ with $t_0\leqslant t \leqslant 2t_0+h$,
$$
\begin{array}{lll}
\displaystyle|(Py)(t)|\leqslant
\eta+(K_1+K_2\eta)\frac{(t_0+h)^\alpha}{\Gamma(\alpha+1)}.
\end{array}
$$
Hence $||P(y)||_{\infty}\leqslant\ell$.

\noindent Step 3. $P$ maps bounded sets into equicontinuous sets on
every compact subset $[t_0-h, b]$ of $[t_0-h, \infty)$.

Let $t_1, t_2 \in [t_0, b]$ , $t_1 < t_2$, and let $B_\eta$ be a
bounded set of $BC$ as in Step 2. Let $y\in B_\eta$. Then we have
$$
\begin{array}{l}
\displaystyle|(Py)(t_2)-(Py)(t_1)| \leqslant |y(t_0)e^{-\beta(t_2-t_0)}-y(t_0)e^{-\beta(t_1-t_0)}|\\
\displaystyle~~~~+\frac{1}{\Gamma(\alpha)}\int\limits_{t_0}^{t_1}
\left|\left((t_2-s)^{\alpha-1}e^{-\beta
(t_2-s)}-(t_1-s)^{\alpha-1}e^{-\beta (t_1-s)}\right)
f(s,y_{s})\right|ds\\
\displaystyle~~~~ +\frac{1}{\Gamma(\alpha)}\int\limits_{t_1}^{t_2}
\left|(t_2-s)^{\alpha-1}e^{-\beta (t_2-s)}f(s,y_{s})\right|ds\\
\displaystyle~~\leqslant |y(t_0)|e^{\beta t_0}|e^{-\beta
t_2}-e^{-\beta
t_1}|+\frac{K_1+K_2\eta}{\Gamma(\alpha+1)}(t_2-t_1)^\alpha\\
\displaystyle~~~~+\frac{K_1+K_2\eta}{\Gamma(\alpha)}\int\limits_{t_0}^{t_1}
\left((t_1-s)^{\alpha-1}e^{-\beta
(t_1-s)}-(t_2-s)^{\alpha-1}e^{-\beta (t_2-s)}\right)ds.
\end{array} \eqno(3.18)
$$
Observing that
$$
\begin{array}{lll}
\displaystyle\frac{K_1+K_2\eta}{\Gamma(\alpha)}\int\limits_{t_0}^{t_1}
\left((t_1-s)^{\alpha-1}e^{-\beta (t_1-s)}-(t_2-s)^{\alpha-1}e^{-\beta (t_1-s)}\right)ds\\
\displaystyle~~\leqslant
\frac{K_1+K_2\eta}{\Gamma(\alpha)}\int\limits_{t_0}^{t_1}
\left((t_1-s)^{\alpha-1}-(t_2-s)^{\alpha-1}\right)ds\\
\displaystyle~~\leqslant\frac{K_1+K_2\eta}{\Gamma(\alpha+1)}
\left((t_1-t_0)^{\alpha}-(t_2-t_0)^{\alpha}+(t_2-t_1)^\alpha\right)\\
\displaystyle~~\leqslant\frac{K_1+K_2\eta}{\Gamma(\alpha+1)}(t_2-t_1)^\alpha,
\end{array} \eqno(3.19)
$$
and from Taylor's theorem, we obtain
$$
\begin{array}{lll}
\displaystyle\frac{K_1+K_2\eta}{\Gamma(\alpha)}\int\limits_{t_0}^{t_1}
\left((t_2-s)^{\alpha-1}e^{-\beta (t_1-s)}-(t_2-s)^{\alpha-1}e^{-\beta (t_2-s)}\right)ds\\
\displaystyle~~\leqslant
\frac{K_1+K_2\eta}{\Gamma(\alpha)}(t_2-t_1)^{\alpha-1}\int\limits_{t_0}^{t_1}
\left(e^{-\beta (t_1-s)}-e^{-\beta (t_2-s)}\right)ds\\
\displaystyle~~\leqslant\frac{K_1+K_2\eta}{\beta\Gamma(\alpha)}(t_2-t_1)^{\alpha-1}
\left(1-e^{-\beta(t_2-t_1)}\right)\\
\displaystyle~~=\frac{K_1+K_2\eta}{\Gamma(\alpha)}\left((t_2-t_1)^\alpha
+\frac{o(t_2-t_1)}{t_2-t_1}(t_2-t_1)^\alpha\right),
\end{array} \eqno(3.20)
$$
where $\lim_{t_2-t_1 \rightarrow 0}\frac{o(t_2-t_1)}{t_2-t_1}=0$. By
(3.18)-(3.20), we can conclude that
$$
\begin{array}{l}
\displaystyle|(Py)(t_2)-(Py)(t_1)| \leqslant \eta e^{\beta t_0}|e^{-\beta t_2}-e^{-\beta t_1}|\\
\displaystyle~~+\frac{2(K_1+K_2\eta)}{\Gamma(\alpha+1)}(t_2-t_1)^\alpha
+\frac{K_1+K_2\eta}{\Gamma(\alpha)}\left((t_2-t_1)^\alpha
+\frac{o(t_2-t_1)}{t_2-t_1}(t_2-t_1)^\alpha\right).
\end{array}
$$
As $t_1 \rightarrow t_2$ the right-hand side of the above inequality
tends to zero. The equicontinuity for the cases $t_1 < t_2 \leqslant
t_0$ and $t_1 \leqslant t_0 \leqslant t_2$ is obvious.

 \noindent Step 4. $P$ maps bounded sets into equiconvergent sets.

 Let  $y\in B_\eta$. Then
$$
\begin{array}{l}
\displaystyle|(Py)(t)|\leqslant|y(t_0)|e^{-\beta (t-t_0)}+
\frac{1}{\Gamma(\alpha)}\int\limits_{t_0}^t
(t-s)^{\alpha-1}e^{-\beta (t-s)}|f(s,y_{s})|ds\\
\displaystyle~~\leqslant \eta e^{-\beta
(t-t_0)}+\frac{1}{\Gamma(\alpha)}\int\limits_{t_0}^t
(t-s)^{\alpha-1}e^{-\beta (t-s)}\left(k_1(s)+k_2(s)\eta\right)ds.
\end{array}
$$
Therefore $(H_4)$ implies that   $|(Py)(t)|$ uniformly (w.r.t $y \in
B(\eta)$) converges to $0$ as $t \rightarrow \infty$. As a
consequence of Steps 1-4, we can conclude that $P : BC \rightarrow
BC$ is continuous and completely continuous.

 \noindent Step 5 (A priori bounds). We now show there exists an open set $U
\subseteq BC$ with $y\neq\lambda P(y)$ for $\lambda \in (0,1)$ and
$y \in
\partial U$.

Let $y\in BC$ and $y =\lambda P(y)$ for some $0 < \lambda < 1$. Then
for each $t \in [t_0,\infty)$ we obtain
$$
y(t) = \lambda\left[y(t_0)e^{-\beta
(t-t_0)}+\frac{1}{\Gamma(\alpha)}\int\limits_{t_0}^t
(t-s)^{\alpha-1}e^{-\beta (t-s)}f(s,y_s)ds\right].
$$
By $(H_3)$, we have that for all $\theta \in [-h, 0]$ and $t
\geqslant t_0+h$,
$$
\begin{array}{lll}
\displaystyle|y(t+\theta)| &\leqslant &
\displaystyle|y(t_0)|e^{-\beta
(t+\theta-t_0)}+\frac{1}{\Gamma(\alpha)}\int\limits_{t_0}^{t+\theta}
(t+\theta-s)^{\alpha-1}e^{-\beta (t+\theta-s)}|f(s,y_s)|ds\\
&\leqslant & \displaystyle|y(t_0)|e^{-\beta
(t+\theta-t_0)}+\frac{1}{\Gamma(\alpha)}\int\limits_{t_0}^{t+\theta}
(t+\theta-s)^{\alpha-1}e^{-\beta
(t+\theta-s)}\left(K_1+K_2\|y_s\|\right)ds,
\end{array}
$$
and thus
$$
\begin{array}{lll}
\displaystyle\|y_t\| &\leqslant & \displaystyle|y(t_0)|e^{-\beta
(t-h-t_0)}+\frac{1}{\Gamma(\alpha)}\int\limits_{t_0}^{t}
(t-s)^{\alpha-1}e^{-\beta (t-h-s)}\left(K_1+K_2\|y_s\|\right)ds.
\end{array}
$$
It follows from the arguments in (3.12)-(3.13), we can conclude that
for each $t \in [t_0, \infty)$,
$$
\displaystyle \frac{1}{\Gamma(\alpha)}\int\limits_{t_0}^{t}
(t-s)^{\alpha-1}e^{-\beta (t-s)}ds \leqslant
\frac{(t_0+h)^{\alpha-1}e^{-\beta (t_0+h)}}{\beta\Gamma(\alpha)}+
\frac{(t_0+h)^{\alpha}}{\Gamma(\alpha+1)}=: R_1.
$$
Hence
$$
\begin{array}{lll}
\displaystyle e^{\beta t}\|y_t\| &\leqslant &
\displaystyle\|\phi\|e^{\beta (h+t_0)}+e^{\beta
h}K_1R_1+\frac{K_2e^{\beta h}}{\Gamma(\alpha)}\int\limits_{t_0}^{t}
(t-s)^{\alpha-1}e^{\beta s}\|y_s\|ds.
\end{array}
$$
Let $R_2=\|\phi\|e^{\beta (h+t_0)}+e^{\beta h}K_1R_1$. Then from
Lemma 2.2, there exists $K$ such that we have for all $t \geqslant
t_0+h$,
$$
\|y_t\| \leqslant R_2+\frac{K K_2 R_2e^{\beta
h}}{\Gamma(\alpha+1)}(t-t_0)^\alpha e^{-\beta t}.
$$
Since $\lim_{t \rightarrow \infty}(t-t_0)^\alpha e^{-\beta t}=0$,
there exists $R_3>0$ such that $$\|y\|_{\infty} \leqslant R_3.$$ Set
$$U=\{y\in BC~:~ ||y||_{\infty}<R_3+1\}.$$
$P : U \rightarrow BC$ is continuous and completely continuous. From
the choice of $U$, there is no $y \in \partial U$ such that $y =
\lambda P(y)$, for $\lambda \in (0,1)$. As a consequence of
Leray-Schauder fixed point theorem, we deduce that $P$ has a fixed
point $y$ in $U$.

\noindent Step 6. Uniform asymptotic stability of solutions.

Let $B \subset C([-h, 0], \mathbb{R})$ be bounded, i.e., there
exists $d \geqslant 0$ such that
$$
\|\psi\|=\sup \limits_{\theta \in [-h, 0]}|\psi(\theta)| \leqslant d
~~~\mbox{ for all } \psi \in B.
$$
From the similar arguments in step 4, we can deduce that there
exists $R_4>0$ such that for all solutions $y(t, t_0, \phi)$ of IVP
(1.1) with initial data $\phi \in B$, we have
$$
\|y\|_{\infty} \leqslant R_4, ~~~\forall \phi \in B.
$$

Now we consider  two solutions $x=x(t)$ and $y=y(t)$ of IVP (1.1)
corresponding to initial values $\psi$ and $\phi$. Note that for all
$t \geqslant t_0$,
$$
\begin{array}{l}
\displaystyle|x(t)-y(t)|\leqslant|x(t_0)-y(t_0)|e^{-\beta (t-t_0)}+
\frac{1}{\Gamma(\alpha)}\int\limits_{t_0}^t
(t-s)^{\alpha-1}e^{-\beta (t-s)}\left(|f(s,x_s)|+|f(s,y_s)|\right)ds\\
~~ \leqslant\displaystyle 2de^{-\beta (t-t_0)}+
\frac{2}{\Gamma(\alpha)}\int\limits_{t_0}^t
(t-s)^{\alpha-1}e^{-\beta (t-s)}\left(k_1(s)+k_2(s)R_4\right)ds.
\end{array} \eqno(3.21)
$$
Then the proof of uniform asymptotic stability of solutions can be
done by making use of $(H_4)$ and (3.21).

The proof of theorem 3.2 is completed. $\Box$

\bigbreak

\noindent{\section{Examples}}

\vskip 5pt \noindent {\bf Example  4.1.} Consider the fractional
functional differential equation
$$
\left\{\begin{array}{l} \displaystyle D^{\frac{1}{2}} \left[y(t)e^{
t}\right]
=\frac{e^{2t}}{8(e^{t}+e^{-t})}\sin^4\left(y\left(t-1\right)\right)+e^t,
\quad t \geqslant 0, \\
y(t)=\phi(t),\quad -1 \leqslant t \leqslant 0,
\end{array}
\right. \eqno(4.1)
$$
where $f(t,
y_t)=\frac{e^{t}}{8(e^{t}+e^{-t})}\sin^4\left(y\left(t-1\right)\right)+1$.
It is clear that condition $(H_0)$ holds. Let $x_t$, $y_t \in C([-1,
0], \mathbb{R})$. Then for all $t \in [0, \infty)$, we have
$$
\begin{array}{lll}
\displaystyle |f(t, x_t)-f(t, y_t)|& = & \displaystyle
\frac{e^{t}}{8(e^{t}+e^{-t})}\left|\sin^4\left(x\left(t-1\right)\right)
-\sin^4\left(y\left(t-1\right)\right)\right|\\
& \leqslant & \displaystyle
\frac{e^{t}}{2(e^{t}+e^{-t})}|x(t-1)-y(t-1)| \leqslant \frac{1}{2}
|x(t-1)-y(t-1)|.
\end{array}
$$
On the other hand, note that $f(t, 0)=1$ for each $t \in [0,
\infty)$ and
$\frac{1}{2}\left(\frac{e^{-1}}{\Gamma(\frac{1}{2})}+\frac{1}{\Gamma(\frac{3}{2})}\right)<1$.
Hence conditions $(H_1)$ and $(3.2)$ hold. By Theorem 3.1, we
conclude that IVP (4.1) has a unique solution in the space $BC([-1,
\infty), \mathbb{R})$, and the solutions of IVP (4.1) are uniformly
asymptotically stable.

\vskip 5pt \noindent {\bf Example  4.2.} Consider the fractional
functional differential equation
$$
\left\{\begin{array}{l} \displaystyle D^{\frac{1}{2}} \left[y(t)e^{
t}\right] =10e^t(t+1)^{-\frac{3}{4}}\frac{y(t-1)}{1+|y(t-1)|},
\quad t \geqslant 0, \\
y(t)=\phi(t),\quad -1 \leqslant t \leqslant 0,
\end{array}
\right. \eqno(4.2)
$$
where $f(t, y_t)=10(t+1)^{-\frac{3}{4}}\frac{y(t-1)}{1+|y(t-1)|}$.
It is easy to see that condition $(H_2)$ holds. Let $y_t \in C([-1,
0], \mathbb{R})$. Then for all $t \in [0, \infty)$, we find that
$$
\displaystyle |f(t, y_t)| =  \displaystyle
\left|10(t+1)^{-\frac{3}{4}}\frac{y(t-1)}{1+|y(t-1)|}\right|\leqslant
10(t+1)^{-\frac{3}{4}}|y(t-1)|,
$$
where  $10(t+1)^{-\frac{3}{4}} \in BC([0, \infty), \mathbb{R}_+)$
with $\sup_{t \geqslant 0}10(t+1)^{-\frac{3}{4}}=10$, and
$$
\frac{1}{\Gamma(\frac{1}{2})}\int\limits_{0}^t (t-s)^{-\frac{1}{2}}
e^{-(t-s)}10(s+1)^{-\frac{3}{4}}ds \leqslant
\frac{10}{\Gamma(\frac{1}{2})}\int\limits_{0}^t (t-s)^{-\frac{1}{2}}
s^{-\frac{3}{4}}ds=\frac{10\Gamma (\frac{1}{4})}{\Gamma
(\frac{3}{4})}t^{-\frac{1}{4}} \rightarrow 0 \mbox{ as } t
\rightarrow \infty.
$$
Thus  conditions $(H_3)$ and $(H_4)$ hold, and the global existence
and  the uniform asymptotic stability of solutions of IVP (4.2) can
be obtained by applying Theorem 3.2.

By using the algorithm given in \cite{Deng}, we numerically simulate
Example 1 with the initial conditions
$\phi(t)=\sin(t),\cos(t),-\cos(t),1.5$, and Example 2 with
$\phi(t)=t,\cos(t),-\cos(t),1.5$, see Figures 1 and 2. From the
numerical results, it can be noted that both of the solutions of
Examples 1 and 2 converge uniformly, and the solutions of Example 1
converge faster than the ones of Example 2. The numerical results
confirm the theoretical analysis.

\begin{figure}
\centering

\includegraphics[width=1.0\textwidth,height=8.8cm]{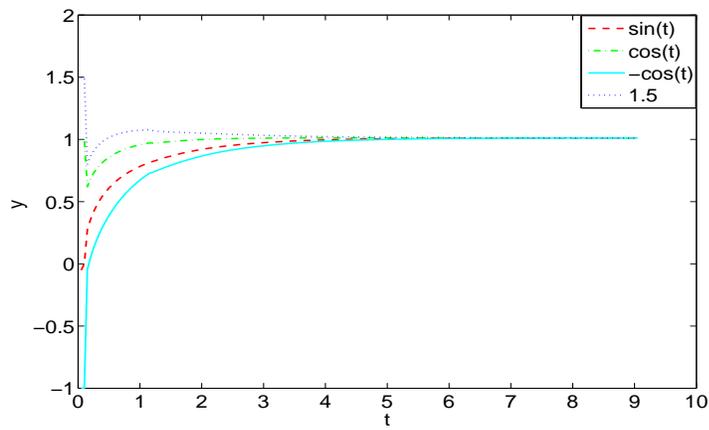}

\caption{The numerical solutions of Example 1 with the initial conditions $\phi(t)=\sin(t),\cos(t),-\cos(t),1.5$, respectively.}%
\end{figure}

\begin{figure}
\centering

\includegraphics[width=1.0\textwidth,height=8.8cm]{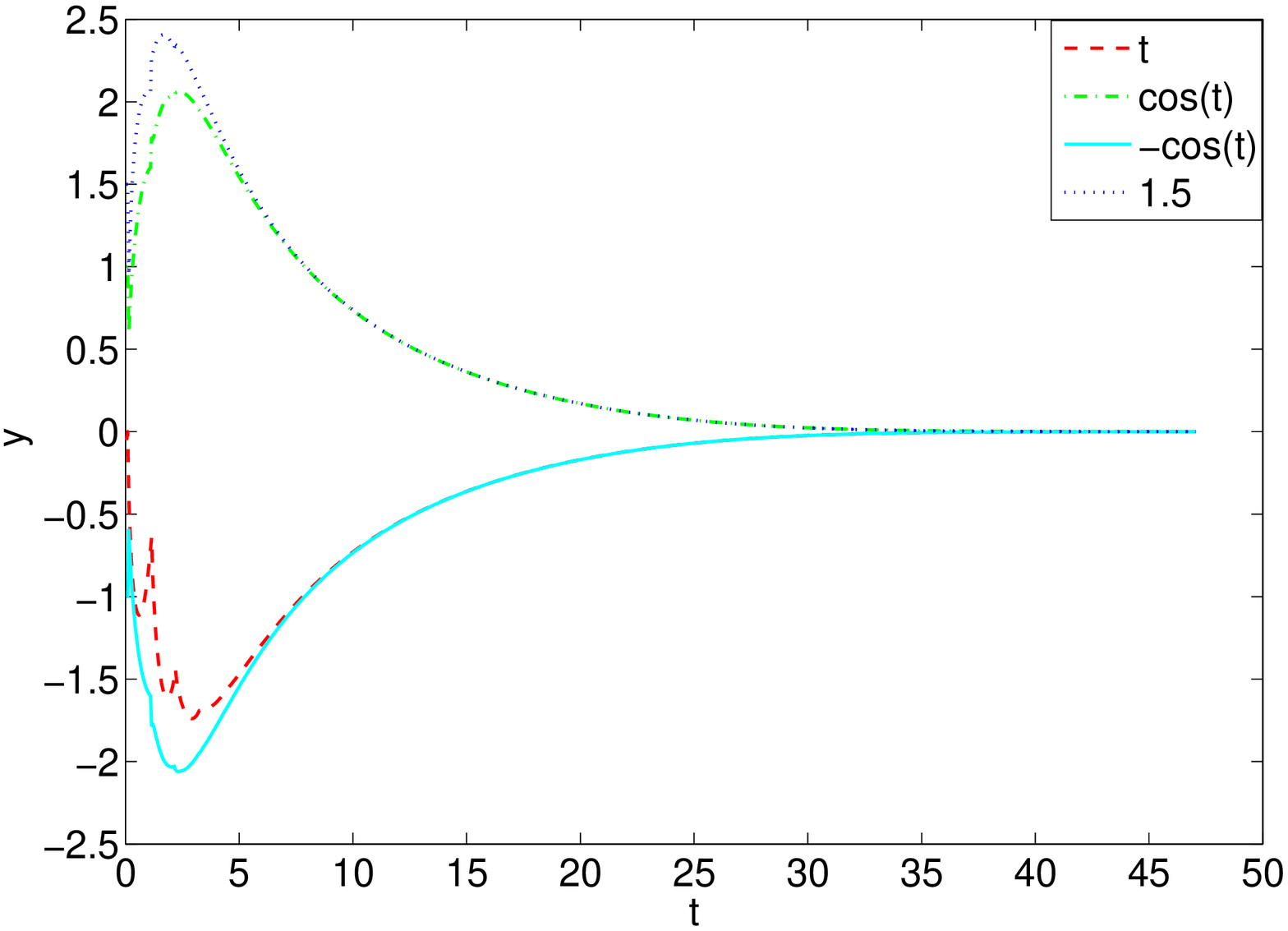}

\caption{The numerical solutions of Example 2 with the initial conditions $\phi(t)=t,\cos(t),-\cos(t),1.5$, respectively.}%
\end{figure}
\bigbreak

\end{document}